\documentclass[11pt]{article}
\usepackage{enumerate}
\usepackage{amssymb,a4wide,latexsym,makeidx,epsfig,fleqn}
\usepackage{amsthm}
\usepackage{amsmath}
\usepackage{enumerate}
\usepackage{graphicx}
\usepackage{subfigure}
\usepackage{float}
\usepackage{caption}
\usepackage[colorlinks, linkcolor=black, anchorcolor=black, citecolor=blue]{hyperref}%
\newtheorem{theorem}{Theorem}[section]

\newtheorem{definition}[theorem]{Definition}
\newtheorem{lemma}[theorem]{Lemma}

\newtheorem{problem}[theorem]{Problem}

\begin{document}
\textwidth 150mm \textheight 225mm
\title{A sharp upper bound on the spectral radius of $\theta(1,3,3)$-free
graphs with given size\thanks{Supported by the National Natural Science Foundation of China (No. 12271439).}}
\author{{Yuxiang Liu$^{a,b}$, Ligong Wang$^{a,b,}$\footnote{Corresponding author.}}\\
{\small $^{a}$  School of Mathematics and Statistics, Northwestern Polytechnical University,}\\{\small  Xi'an, Shanxi 710129, P.R. China.}\\ {\small $^{b}$ Xi'an-Budapest Joint Research Center for Combinatorics, Northwestern Polytechnical University,}\\{\small  Xi'an, Shanxi 710129, P.R. China.}\\
{\small E-mail: yxliumath@163.com, lgwangmath@163.com}}
\date{}
\maketitle
\begin{center}
\begin{minipage}{135mm}
\vskip 0.3cm
\begin{center}
{\small {\bf Abstract}}
\end{center}
{\small  A graph $G$ is $F$-free if $G$ does not contain $F$ as a subgraph. Let $\rho(G)$ be the spectral radius of a graph $G$. Let $\theta(1,p,q)$ denote the theta graph, which is obtained by connecting two distinct vertices with three internally disjoint paths with lengths $1, p, q$, where $p\leq q$. Let $S_{n,k}$ denote the graph obtained by joining every vertex of $K_{k}$ to $n-k$ isolated vertices and $S_{n,k}^{-}$ denote the graph obtained from $S_{n,k}$ by deleting an edge incident to a vertex of degree $k$, respectively. In this paper, we show that if $\rho(G)\geq\rho(S_{\frac{m+4}{2},2}^{-})$ for a graph $G$ with even size $m\geq 92$, then $G$ contains a $\theta(1,3,3)$ unless $G\cong S_{\frac{m+4}{2},2}^{-}$.
}
\end{minipage}
\end{center}
\section{Introduction}
Throughout this paper, we consider all graphs are always undirected and simple. We follow the traditional notation and terminology \cite{BoMu1}. Let $G$ be a graph of order $n$ with vertex set $V(G)=\{v_{1}, v_{2},\ldots, v_{n}\}$ and size $m$ with edge set $E(G)$. For a vertex $u\in V(G)$, let $N_{G}(u)$ be the neighborhood set of a vertex $u$, $N_{G}[u]=N_{G}(u)\cup \{u\}$ and $N_{G}^{2}(u)$ be the set of vertices of distance two to $u$ in $G$. In particular, $N_{S}(v)=N(v)\cap S$ and $d_{S}(v)=|N_{S}(v)|$ for a subset $S\subseteq V(G)$. Let $d_{G}(u)=|N_{G}(u)|$ be the degree of a vertex $u$. For the sake of simplicity, we omit all the subscripts if $G$ is clear from the context, for example, $N(u)$, $N[u]$, $N^{2}(u)$ and $d(u)$. For a graph $G$ and a subset $S\subseteq V(G)$, let $G[S]$ be the subgraph of $G$ induced by $S$. For two vertex subsets $S$ and $T$ of $V(G)$ (where $S\cap T$ may not be empty), let $e(S,T)$ denote the number of edges with one endpoint in $S$ and the other in $T$. $e(S,S)$ is simplified by $e(S)$. Given two vertex-disjoint graphs $G_{1}$ and $G_{2}$, we denote by $G_{1}\cup G_{2}$ the disjoint union of the two graphs, and by $G_{1}\vee G_{2}$ the joint graph obtained from $G_{1}\cup G_{2}$ by joining each vertex of $G_{1}$ with each vertex of $G_{2}$. The adjacency matrix of a graph $G$ is an $n\times n$ matrix $A(G)$ whose $(i,j)$-entry is $1$ if $v_{i}$ is adjacent to $v_{j}$ and $0$ otherwise. The spectral radius $\rho(G)$ of a graph $G$ is the largest eigenvalue of its adjacency matrix $A(G)$.

Let $P_{n}, C_{n}, K_{1,n-1}$ and $K_{a,b}$ be the path of order $n$, the cycle of order $n$, the star graph of order $n$ and the complete bipartite graph with two parts of sizes $a, b$, respectively. Let $S_{n}^{k}$ be the graph obtained from $K_{1,n-1}$ by adding $k$ disjoint edges within its independent sets. Let $S_{n,k}$ be the graph obtained by joining every vertex of $K_{k}$ to $n-k$ isolated vertices. Let $S_{n,k}^{-}$ be the graph obtained from $S_{n,k}$ by deleting an edge incident to a vertex of degree $k$. Let $\theta(1,p,q)$ denote the theta graph, which is obtained by connecting two distinct vertices with three internally disjoint paths with lengths $1, p, q$, where $q\leq p$.

Given a graph $F$, a graph $G$ is $F$-free if it does not contain $F$ as a subgraph. Let $\mathcal{G}(m, F)$ denote the family of $F$-free graphs with $m$ edges and without isolated vertices. The classic Tur\'{a}n type problem asks what is the maximum number of edges in an $F$-free graph of order $n$. In spectral graph theory, Nikiforov \cite{Ni2} posed a spectral Tur\'{a}n type problem which asks to determine the maximum spectral radius of an $F$-free graphs of $n$ vertices, which is known as the Brualdi-Solheid-Tur\'{a}n type problem. In the past decades, this problem has received much attention. For more details, we suggest the reader to see surveys \cite {ChZ,FSi,LiLF,Ni3}, and references therein. In addition, Brualdi and Hoffman raised another spectral Tur\'{a}n type problem: What is the maximal spectral radius of an $F$-free graph with given size $m$? This problem is called the Brualdi-Hoffman-Tur\'{a}n type problem. This problem has been studied for various families of graphs. For example, $K_{3}$ \cite{Nos}, $K_{r+1}$ \cite{Ni5,Ni6}, $K_{2,r+1}$ \cite{ZhLS}, $F_{2k+2}$ \cite{LiZZ} (where $F_{k}=K_{1}\vee P_{k-1}$), $F_{k,3}$ \cite{LiZZ} (where $F_{k,3}$ is the friendship graph obtained from $k$ triangles by sharing a common vertex).

For theta graphs, Sun, Li and Wei \cite{SLW} characterized the extremal graph with maximum spectral radius of $\theta(1,2,3)$-free and $\theta(1,2,4)$-free graphs with odd size. Fang and You \cite{FaY} characterized the extremal graph with maximum spectral radius of $\theta(1,2,3)$-free graphs with even size. Liu and Wang \cite{LiuW} characterized the extremal graph with maximum spectral radius of $\theta(1,2,4)$-free graphs with even size. Lu, Lu and Li \cite{LuLL} characterized the extremal graph with maximum spectral radius of $\theta(1,2,5)$-free with given size. Li, Zhai and Shu \cite{LiZS} characterized the extremal graph with maximum spectral radius of $\theta(1,2,2k-1)$-free or $\theta(1,2,2k)$-free with given size.

Recently, Li, Zhao and Zou \cite{LiZZ} characterized the extremal graph with maximum spectral radius of $\theta(1,p,q)$-free with size $m$ for $q\geq p\geq 3$ and $p+q\geq7$.

\noindent\begin{theorem}\label{th:ch-1.1.}{\rm(}$\cite{LiZZ}${\rm)}
Let $k\geq3$ and $m\geq\frac{9}{4}k^{6}+6k^{5}+46k^{4}+56k^{3}+196k^{2}$. If $G\in \mathcal{G}(m, \theta(1,p,q))\cup  \mathcal{G}(m, \theta(1,r,s))$ with $q\geq p\geq 3$, $s\geq r\geq3$, $p+q=2k+1$ and $r+s=2k+2$, then $$\rho(G)\leq\frac{k-1+\sqrt{4m-k^{2}+1}}{2},$$ and equality holds if and only if $G\cong K_{k}\vee (\frac{m}{k}-\frac{k-1}{2})K_{1}$.
\end{theorem}

At the same time, they \cite{LiZZ} proposed the following problem.

\noindent\begin{problem}\label{pr:ch-1.2.}{\rm(}$\cite{LiZZ}${\rm)}
How can we characterize the graphs among $\mathcal{G}(m,\theta(1,3,3))$ having the largest spectral radius?
\end{problem}

\noindent\begin{theorem}\label{th:ch-1.3.}{\rm(}$\cite{GaoL}${\rm)}
Let $G\in \mathcal{G}(m, \theta(1,3,3))$ be a graph of size $m\geq43$. Then $$\rho(G)\leq \frac{1+\sqrt{4m-3}}{2},$$ and equality holds if and only if $G\cong S_{\frac{m+3}{2},2}$.
\end{theorem}

However, for a $\theta(1,3,3)$-free graph $G$ with size $m$, the bound $\rho(G)\leq \frac{1+\sqrt{4m-3}}{2}$ is sharp only for odd $m$. Motivated by this, we want to obtain a sharp upper bound of $\rho(G)$ for even $m$. Our result is presented as follows.

\noindent\begin{theorem}\label{th:ch-1.4.} Let $G\in \mathcal{G}(m, \theta(1,3,3))$ be a graph of even size $m\geq92$, then $\rho(G)\leq \rho(S_{\frac{m+4}{2},2}^{-})$, and equality holds if and only if $G\cong S_{\frac{m+4}{2},2}^{-}$.
\end{theorem}

\section{Preliminary}
In this section, we introduce some lemmas which are used to prove our result.

\noindent\begin{lemma}\label{le:ch-2.1.} {\rm(}$\cite{ZhWF1}${\rm)} Let $u, v$ be two distinct vertices of the connected graph $G$, $\{v_{i}| i=1,2, \ldots, s\}\subseteq N(v)\setminus N(u)$, and $X=(x_{1}, x_{2}, \ldots, x_{n})^{T}$ be the Perron vector of $G$. Let $G^{\prime}=G-\sum_{i=1}^{s}v_{i}v+\sum_{i=1}^{s}v_{i}u$. If $x_{u}\geq x_{v}$, then $\rho(G)<\rho(G^{\prime})$.
\end{lemma}

\noindent\begin{definition}\label{de:ch-2.2.}{\rm(}$\cite{CvRS}${\rm)} Given a graph $G$, the vertex partition $\Pi$: $V(G)=V_{1}\cup V_{2} \cup \ldots \cup V_{k}$ is said to be an equitable partition if, for each $u\in V_{i}$, $|V_{j}\cap N(u)|=b_{ij}$ is a constant depending only on $i,j$ ($1\leq i,j\leq k$). The matrix $B_{\Pi}=(b_{ij})$ is called the quotient matrix of $G$ with respect to $\Pi$.
\end{definition}

\noindent\begin{lemma}\label{le:ch-2.3.}{\rm(}$\cite{CvRS}${\rm)} Let $\Pi$: $V(G)=V_{1}\cup V_{2} \ldots \cup V_{k}$ be an equitable partition of a graph $G$ with quotient matrix $B_{\Pi}$. Then $det(xI-B_{\Pi}) \mid det(xI-A(G))$. Furthermore, the largest eigenvalue of $B_{\Pi}$ is just the spectral radius of $G$.
\end{lemma}

\noindent\begin{lemma}\label{le:ch-2.4.}{\rm(}$\cite{ZhLS}${\rm)} Let $G^{\ast}$ be the extremal graph with the maximum spectral radius in $\mathcal{G}(m,F)$. Let $X=(x_{1}, x_{2}, \ldots, x_{n})^{T}$ be the Perron vector of the graph $G^{\ast}$. If $F$ is a $2$-connected graph and $x_{u^{\ast}}=max\{x_{v}\mid v\in V(G^{\ast})\}$, then the following statements hold.

{\rm(}$i${\rm)} $G^{\ast}$ is connected.

{\rm(}$ii${\rm)} There exists no cut vertex in $V(G^{\ast})\setminus\{u^{\ast}\}$, and hence $d(u)\geq2$ for any $u\in V(G^{\ast})\setminus N[u^{\ast}]$.
\end{lemma}

\noindent\begin{lemma}\label{le:ch-2.5.}{\rm(}$\cite{BrH}${\rm)} Let $G$ be a bipartite graph of size $m$. Then $\rho(G)\leq\sqrt{m}$, and equality holds if and only if $G$ is a disjoint union of a complete bipartite graph and isolated vertices.
\end{lemma}

\noindent\begin{lemma}\label{le:ch-2.6.}{\rm(}$\cite{MiLH}${\rm)} $\rho(S_{\frac{m+4}{2},2}^{-})>\frac{1+\sqrt{4m-5}}{2}$ for $m\geq6$.

\end{lemma}

\noindent\begin{lemma}\label{le:ch-2.7.}{\rm(}$\cite{MiLH}${\rm)} Let $X=(x_{1}, x_{2}, \ldots, x_{n})^{T}$ be the Perron vector of a connected graph $G$ of size $m$ and let $x_{u^{\star}}=max\{x_{v}| v\in V(G)\}$ and $W=V(G)\backslash N[u^{\star}]$. If $\rho(G)>\frac{1+\sqrt{4m-5}}{2}$ and there exists a vertex $v$ of $G$ such that $x_{v}< (1-\beta)x_{u^{\star}}$, where $0<\beta<1$, then
$$e(W)<e(N(u^{\star}))-|N(u^{\star})\setminus N_{0}(u^{\star})|+\frac{3}{2}-\beta d_{N(u^{\star})}(v),$$  $\mbox{for }  v\in N^{2}(u^{\star})$.
\end{lemma}

\section{Proof of Theorem 1.4.}
Let $G^{\ast}$ be the extremal graph with maximum spectral radius in $\mathcal{G}(m,F)$. Let $\rho=\rho(G^{\ast})$ and $X^{\ast}=(x_{1}, x_{2}, \ldots, x_{n})^{T}$ be the Perron vector of $G^{\ast}$ with coordinate $x_{v}$ corresponding to the vertex $v\in V(G^{\ast})$. A vertex $u^{\ast}$ in $G^{\ast}$ is said to be an extremal vertex if $x_{u^{\ast}}=max\{x_{v}\mid v\in V(G^{\ast})\}$. Let $W=V(G^{\ast})\backslash N[u^{\ast}]$ and $N_{+}(u^{\ast})=N(u^{\ast})\setminus N_{0}(u^{\ast})$, where $N_{0}(u^{\ast})$ denotes the set of isolated vertices of $G^{\ast}[N(u^{\ast})]$. Let $W_{H}=\cup_{u\in V(H)}N_{W}(u)$ for any component $H$ of $G^{\ast}[N(u^{\ast})]$.

\noindent\begin{lemma}\label{le:ch-3.1.} Let $G^{\ast}\in \mathcal{G}(m, \theta(1,3,3))$. Then $G^{\ast}[N(u^{\ast})]$ is $P_{5}$-free, that is, each component $H$ of $G^{\ast}[N(u^{\ast})]$ is one of the following.

{\rm(}$i${\rm)} a graph with $C_{4}$ as its spanning subgraph, that is, $C_{4}$, $\theta(1,2,2)$ or $K_{4}$;

{\rm(}$ii${\rm)} a copy of $S_{r+1}^{1}$ for $r\geq2$, where $S_{3}^{1}$ is a triangle for $r=2$;

{\rm(}$iii${\rm)} a double star $D_{a,b}$ for $a,b\geq1$, which is obtained from two stars $K_{1,a}$ and $K_{1,b}$
by joining a new edge between their centers;

{\rm(}$iv${\rm)} a star $K_{1,r}$ for $r\geq0$, where $K_{1,0}$ is a singleton component.

\end{lemma}

\noindent\begin{lemma}\label{le:ch-3.2.}{\rm(}$\cite{GaoL}${\rm)}
For any non-trivial component $H$ in $G^{\ast}[N(u^{\ast})]$, if $H$ contains a cycle of length four, then $N_{W}(u)\cap N_{W}(v)=\emptyset$ for any two vertices $u$ and $v$ in the cycle of length four.
\end{lemma}

Since $A(G^{\ast})X=\rho X$, we have

$$\rho x_{u^{\ast}}=\sum_{u\in N_{+}(u^{\ast})}x_u+\sum_{u\in N_{0}(u^{\ast})}x_{u}.$$

Furthermore, we have
$$
\begin{aligned}
\rho^{2}x_{u^{\ast}}&=\sum_{u\in N_{+}(u^{\ast})}\rho x_u+\sum_{u\in N_{0}(u^{\ast})}\rho x_{u}
\\&=|N_{+}(u^{\ast})|x_{u^{\ast}}+\sum_{u\in N_{+}(u^{\ast})}d_{N(u^{\ast})}(u)x_{u}+\sum_{w\in N_{W}(u)}d_{W}(u)x_w+|N_{0}(u^{\ast})|x_{u^{\ast}}
\\&=d(u^{\ast})x_{u^{\ast}}+\sum_{u\in N_{+}(u^{\ast})}d_{N(u^{\ast})}(u)x_{u}+\sum_{w\in N^{2}(u^{\ast})}d_{N(u^{\ast})}(w)x_w.
\end{aligned}
$$
Therefore,
\begin{equation}\label{eq:ch-1}
\begin{aligned}
(\rho^{2}-\rho)x_{u^{\ast}}&=d(u^{\ast})x_{u^{\ast}}+\sum_{v\in N_{+}(u^{\ast})}(d_{N(u^{\ast})}(v)-1)x_{v}+\sum_{w\in N^{2}(u^{\ast})}d_{N(u^{\ast})}(w)x_w-\sum_{v\in N_{0}(u^{\ast})}x_{v}
\\&\leq |N(u^{\ast})|x_{u^{\ast}}+\sum_{u\in N_{+}(u^{\ast})}(d_{N(u^{\ast})}(u)-1)x_{u}+e(N(u^{\ast}),W)x_{u^{\ast}}-\sum_{u\in N_{0}(u^{\ast})}x_{u}.
\end{aligned}
\end{equation}

Note that $S_{\frac{m+4}{2},2}^{-}$ is $\theta(1,3,3)$-free, we have $\rho\geq\rho(S_{\frac{m+4}{2},2}^{-})> \frac{1+\sqrt{4m-5}}{2}>10$ for $m\geq92$. By Lemma \ref{le:ch-2.6.}, we have

\begin{equation}\label{eq:ch-2}
(\rho^{2}-\rho)x_{u^{\ast}}>(m-\frac{3}{2})x_{u^{\ast}}=(|N(u^{\ast})|+e(N_{+}(u^{\ast}))+e(N(u^{\ast}),W)+e(W)-\frac{3}{2})x_{u^{\ast}}. \end{equation}

Combining with \eqref{eq:ch-1} and \eqref{eq:ch-2}, we get
$$\sum_{u\in N_{+}(u^{\ast})}(d_{N(u^{\ast})}(u)-1)x_{u}>\left(e(N_{+}(u^{\ast}))+e(W)+\sum_{u\in N_{0}(u^{\ast})}\frac{x_{u}}{x_{u^{\ast}}}-\frac{3}{2}\right)x_{u^{\ast}}.$$

Let $\mathcal{H}$ be the set of all non-trivial components in $G^{\ast}[N(u^{\ast})]$. For each non-trivial connected component $H$ of $\mathcal{H}$, we denote $\zeta(H)=\sum_{u\in V(H)}(d_{H}(u)-1)x_{u}$. Obviously,
\begin{equation}\label{eq:ch-3}
\sum_{H\in \mathcal{H}}\zeta(H)>\left(e(N_{+}(u^{\ast}))+e(W)+\sum_{u\in N_{0}(u^{\ast})}\frac{x_{u}}{x_{u^{\ast}}}-\frac{3}{2}\right)x_{u^{\ast}}.
\end{equation}

\noindent\begin{lemma}\label{le:ch-3.3.} $G^{\ast}[N(u^{\ast})]$ contains no any cycle of length four.
\end{lemma}
\noindent{\bf Proof.}  Let $\mathcal{H^{\prime}}$ be the family of components of $G^{\ast}[N(u^{\ast})]$ each of which contains $C_{4}$ as a spanning subgraph and $\mathcal{H}\setminus\mathcal{H^{\prime}}$ be the family of other non-trivial components of $G^{\ast}[N(u^{\ast})]$ each of which contains no $C_{4}$ as a spanning subgraph. By Lemma \ref{le:ch-3.1.} {\rm(}$ii${\rm)}-{\rm(}$iv${\rm)}, for each $H\in\mathcal{H}\setminus\mathcal{H^{\prime}}$, we have
$$\zeta(H)=\sum_{v\in V(H)}(d_{H}(v)-1)x_{v}\leq(2e(H)-|V(H)|)x_{u^{\ast}}\leq e(H)x_{u^{\ast}}.$$

Next we show that $$\zeta(H)<(e(H)-\frac{3}{2})x_{u^{\ast}}+\frac{2\sum_{w\in W_{H}}x_{w}}{\rho-3}$$ for each $H\in\mathcal{H^{\prime}}$. Let $H^{\ast}\in \mathcal{H^{\prime}}$ with $V(H^{\ast})=\{u_{1},u_{2},u_{3},u_{4}\}$ and the cycle of length four be $u_{1}u_{2}u_{3}u_{4}u_{1}$.

First, we consider $W_{H^{\ast}}=\emptyset$. Assume that $x_{u_{1}}=max\{x_{u_{i}}|1\leq i\leq4\}$. Then
$$\rho x_{u_{1}}=\sum_{u\in N(u_{1})}x_{u}\leq x_{u^{\ast}}+x_{u_{2}}+x_{u_{3}}+x_{u_{4}}\leq x_{u^{\ast}}+3x_{u_{1}}.$$
Hence, $x_{u_{1}}\leq \frac{1}{\rho-3}x_{u^{\ast}}<\frac{1}{4}x_{u^{\ast}}$ for $\rho\geq\rho(S_{\frac{m+4}{2},2}^{-})> \frac{1+\sqrt{4m-5}}{2}>10$ due to $m\geq92$. Furthermore, $$\zeta(H^{\ast})\leq(2e(H^{\ast})-|V(H^{\ast})|)x_{u_{1}}<(e(H^{\ast})-\frac{3}{2})x_{u^{\ast}}+\frac{2\sum_{w\in W_{H^{\ast}}}x_{w}}{\rho-3},$$ as desired.

In the following, we assume that $W_{H^{\ast}}\neq\emptyset$, we consider the following two cases.

{\bf Case 1.} All vertices in $W_{H^{\ast}}$ have a unique common neighbor in $V(H^{\ast})$.

Without loss of generality, let $u_{1}\in V(H^{\ast})$ be the unique common neighbor. Thus, $N_{W}(u_{i})=\emptyset$ for $i\in\{2,3,4\}$. Assume that $x_{u_{2}}=max\{x_{u_{i}}|2\leq i\leq4\}$, we have
$$\rho x_{u_{2}}\leq x_{u^{\ast}}+x_{u_{1}}+x_{u_{3}}+x_{u_{4}}\leq 2x_{u_{2}}+2x_{u^{\ast}}.$$
Thus, $x_{u_{2}}\leq \frac{2}{\rho-2}x_{u^{\ast}}<\frac{1}{4}x_{u^{\ast}}$ due to  $\rho>10$. Thus,
$$
\begin{aligned}
\zeta(H^{\ast})&=\sum_{u\in V(H^{\ast})}(d_{H^{\ast}}(u)-1)x_{u}
\\&\leq(d_{H^{\ast}}(u_{1})-1)x_{u_{1}}+(2e(H^{\ast})-d_{H^{\ast}}(u_{1})-3)x_{u_{2}}
\\&< \left(d_{H^{\ast}}(u_{1})-1+\frac{1}{2}e(H^{\ast})-\frac{1}{4}d_{H^{\ast}}(u_{1})-\frac{3}{4}\right)x_{u^{\ast}}
\\&\leq (\frac{3}{4}d_{H^{\ast}}(u_{1})+\frac{1}{2}e(H^{\ast})-\frac{7}{4})x_{u^{\ast}}
\\&\leq(\frac{1}{2}e(H^{\ast})+\frac{1}{2})x_{u^{\ast}}
\\&<(e(H^{\ast})-\frac{3}{2})x_{u^{\ast}}+\frac{2\sum_{w\in W_{H^{\ast}}}x_{w}}{\rho-3}
\end{aligned}
$$
due to $d_{H^{\ast}}(u_{1})\leq 3$, as desired.

{\bf Case 2.} There are at least two distinct vertices of $W_{H^{\ast}}$ such that they have distinct neighbors in $V(H^{\ast})$. Since
\[
\begin{cases}
\rho x_{u_{1}}\leq x_{u^{\ast}}+x_{u_{2}}+x_{u_{3}}+x_{u_{4}}+\sum_{w\in N_{W_{H^{\ast}}}(u_{1})}x_{w},\\
\rho x_{u_{2}}\leq x_{u^{\ast}}+x_{u_{1}}+x_{u_{3}}+x_{u_{4}}+\sum_{w\in N_{W_{H^{\ast}}}(u_{2})}x_{w},\\
\rho x_{u_{3}}\leq x_{u^{\ast}}+x_{u_{1}}+x_{u_{2}}+x_{u_{4}}+\sum_{w\in N_{W_{H^{\ast}}}(u_{3})}x_{w},\\
\rho x_{u_{4}}\leq x_{u^{\ast}}+x_{u_{1}}+x_{u_{2}}+x_{u_{3}}+\sum_{w\in N_{W_{H^{\ast}}}(u_{4})}x_{w},
\end{cases}
\]
we obtain $$\rho(x_{u_{1}}+x_{u_{2}}+x_{u_{3}}+x_{u_{4}})\leq 3(x_{u_{1}}+x_{u_{2}}+x_{u_{3}}+x_{u_{4}})+4x_{u^{\ast}}+\sum_{i=1}^{4}\sum_{w\in N_{W_{H^{\ast}}}(u_{i})}x_{w}.$$
By Lemma \ref{le:ch-3.2.}, we get $N_{W_{H^{\ast}}}(u_{i})\cap N_{W_{H^{\ast}}}(u_{j})=\emptyset$ for arbitrary two distinct vertices $u_{i}, u_{j}\in V(H^{\ast})$. Thus, $\sum_{w\in N_{W_{H^{\ast}}}}x_{w}=\sum_{w\in N_{W}(V(H^{\ast}))}x_{w}=\sum_{i=1}^{4}\sum_{w\in N_{W_{H^{\ast}}}(u_{i})}x_{w}$. Combining $\rho>10$, we obtain
$$
\begin{aligned}
x_{u_{1}}+x_{u_{2}}+x_{u_{3}}+x_{u_{4}}&\leq \frac{4x_{u^{\ast}}}{\rho-3}+\frac{\sum_{w\in W_{H^{\ast}}}x_{w}}{\rho-3}
\\&< \frac{4}{7}x_{u^{\ast}}+\frac{\sum_{w\in W_{H^{\ast}}}x_{w}}{\rho-3}.
\end{aligned}
$$
Hence, by the definition of $\zeta(H^{\ast})$,
$$
\begin{aligned}
\zeta(H^{\ast})&\leq2(x_{u_{1}}+x_{u_{2}}+x_{u_{3}}+x_{u_{4}})
\\&<\frac{8}{7}x_{u^{\ast}}+\frac{2\sum_{w\in W_{H^{\ast}}}x_{w}}{\rho-3}
\\&<(e(H^{\ast})-\frac{3}{2})x_{u^{\ast}}+\frac{2\sum_{w\in W_{H^{\ast}}}x_{w}}{\rho-3}.
\end{aligned}
$$
Thus, we obtain $\zeta(H)<(e(H)-\frac{3}{2})x_{u^{\ast}}+\frac{2\sum_{w\in W_{H}}x_{w}}{\rho-3}$ for each $H\in \mathcal{H^{\prime}}$. Recall that $\zeta(H)\leq e(H)x_{u^{\ast}}$ for each $H\in \mathcal{H}\setminus\mathcal{H^{\prime}}$. Furthermore,
$$
\begin{aligned}
\sum_{H\in \mathcal{H}}\zeta(H)&=\sum_{H\in \mathcal{H^{\prime}}}\zeta(H)+\sum_{H\in \mathcal{H}\setminus\mathcal{H^{\prime}}}\zeta(H)
\\&<\sum_{H\in \mathcal{H^{\prime}}}(e(H)-\frac{3}{2})x_{u^{\ast}}+\sum_{H\in \mathcal{H^{\prime}}}\frac{2\sum_{w\in W_{H}}x_{w}}{\rho-3}+\sum_{H\in \mathcal{H}\setminus\mathcal{H^{\prime}}}e(H)x_{u^{\ast}}
\\&=e(N_{+}(u^{\ast}))x_{u^{\ast}}-\frac{3}{2}\sum_{H\in \mathcal{H^{\prime}}}x_{u^{\ast}}+\sum_{H\in \mathcal{H^{\prime}}}\frac{2\sum_{w\in W_{H}}x_{w}}{\rho-3}.
\end{aligned}
$$
For any $H\in \mathcal{H^{\prime}}$ satisfying $W_{H}=\emptyset$, we have $\sum_{w\in W_{H}}x_{w}=0$. For any $H\in \mathcal{H^{\prime}}$ satisfying $W_{H}\neq\emptyset$ and any $w\in W_{H}$, since $G^{\ast}$ is $\theta(1,3,3)$-free, we get $W_{H}\cap W_{N(u^{\ast})\setminus H}=\emptyset$. Then $d_{N(u^{\ast})\setminus H}(w)=0$. By Lemma \ref{le:ch-3.2.}, we have $d_{H}(w)=1$. Furthermore, $d_{N(u^{\ast})}(w)=1$ and $d_{W}(w)\geq1$ from Lemma \ref{le:ch-2.4.}. Thus, $\sum_{H\in \mathcal{H^{\prime}}}\sum_{w\in W_{H}}x_{w}\leq
\sum_{H\in \mathcal{H^{\prime}}}\sum_{w\in W_{H}}d_{W}(w)x_{w}\leq
\sum_{H\in \mathcal{H^{\prime}}}\sum_{w\in W_{H}}d_{W}(w)x_{u^{\ast}}
\leq 2e(W)x_{u^{\ast}}$. Note that $\rho>10$. Thus,
$$
\begin{aligned}
\sum_{H\in \mathcal{H}}\zeta(H)&<e(N_{+}(u^{\ast}))x_{u^{\ast}}-\frac{3}{2}\sum_{H\in \mathcal{H^{\prime}}}x_{u^{\ast}}+\sum_{H\in \mathcal{H^{\prime}}}\frac{2\sum_{w\in W_{H}}x_{w}}{\rho-3}
\\&<(e(N_{+}(u^{\ast}))+\frac{4}{7}e(W)-\sum_{H\in \mathcal{H^{\prime}}}\frac{3}{2})x_{u^{\ast}}
\\&<(e(N_{+}(u^{\ast}))+e(W)-\sum_{H\in \mathcal{H^{\prime}}}\frac{3}{2})x_{u^{\ast}},
\end{aligned}
$$
which contradicts \eqref{eq:ch-3}. Furthermore, $G^{\ast}[N(u^{\ast})]$ contains no $C_{4}$. This completes the proof of Lemma \ref{le:ch-3.3.}. $\qedsymbol$

By Lemma \ref{le:ch-3.1.}, we know that each non-trivial component of $G^{\ast}[N(u^{\ast})]$ is either a tree or a unicyclic graph $S_{r+1}^{1}$ with $r\geq2$. Let $c$ be the number of non-trivial tree components of $G^{\ast}[N(u^{\ast})]$. Then

$$\sum_{H\in \mathcal{H}}\zeta(H)\leq\sum_{H\in \mathcal{H}}\sum_{u\in V(H)}(d_{H}(u)-1)x_{u^{\ast}}=\sum_{H\in \mathcal{H}}(2e(H)-|H|)x_{u^{\ast}}=(e(N_{+}(u^{\ast}))-c)x_{u^{\ast}}.$$
Combining \eqref{eq:ch-1}, we get

\begin{equation}\label{eq:ch-4}
e(W)<\frac{3}{2}-c-\sum_{u\in N_{0}(u^{\ast})}\frac{x_{u}}{x_{u^{\ast}}}.
\end{equation}

Thus, $e(W)\leq1$ and $c\leq1$. In addition, if $e(W)=1$, then $c=0$ and $\sum_{u\in N_{0}(u^{\ast})}\frac{x_{u}}{x_{u^{\ast}}}<\frac{1}{2}$.

\noindent\begin{lemma}\label{le:ch-3.4.} $e(W)=0$.
\end{lemma}
\noindent{\bf Proof.} Suppose on the contrary that $e(W)=1$. In this case, we have $c=0$ and $\sum_{u\in N_{0}(u^{\ast})}\frac{x_{u}}{x_{u^{\ast}}}$ $<\frac{1}{2}$. It follows that each component of $G^{\ast}[N(u^{\ast})]$ is isomorphic to a unicyclic graph $S_{r+1}^{1}$ with $r\geq2$. That is, each component of $G^{\ast}[N(u^{\ast})]$ contains a triangle. Let $H^{\ast}$ be a component of $G^{\ast}[N(u^{\ast})]$ and $u_{1}u_{2}u_{3}$ is a $C_{3}$ of $H^{\ast}$. Let $w_{1}w_{2}$ be the unique edge of $e(W)$. By Lemma \ref{le:ch-2.4.}, we obtain that $d_{N(u^{\ast})}(w_{i})\geq 1$ for each $i\in \{1,2\}$. If $H^{\ast}\cong S_{3}^{1}$, then we obtain $d_{S_{3}^{1}}(w)\leq3$. If $H^{\ast}\cong S_{r+1}^{1}$ with $r\geq3$, then let $d_{H^{\ast}}(u_{2})=d_{H^{\ast}}(u_{3})=2$ and $u_{4},u_{5},\cdots, u_{r+1}$ be the neighbors of $u_{1}$. For $w\in W_{H^{\ast}}$, we claim $d_{H^{\ast}}(w)\leq1$. Otherwise, we consider the following five cases in the sense of symmetry.
If $\{u_{1},u_{2}\}\subseteq N_{H^{\ast}}(w)$, then $u_{1}u_{3}$,$u_{1}wu_{2}u_{3}$,$u_{1}u_{4}u^{\ast}u_{3}$ are three internally disjoint paths of lengths $1,3,3$ between $u_{1}$ and $u_{3}$, a contradiction. If $\{u_{2},u_{3}\}\subseteq N_{H^{\ast}}(w)$, then $u_{1}u_{3}$,$u_{1}u_{2}wu_{3}$,$u_{1}u_{4}u^{\ast}u_{3}$ are three internally disjoint paths of lengths $1,3,3$ between $u_{1}$ and $u_{3}$, a contradiction. If $\{u_{1}, u_{4}\}\subseteq N_{H^{\ast}}(w)$, then $u_{1}u^{\ast}$,$u_{1}u_{2}u_{3}u^{\ast}$,$u_{1}wu_{4}u^{\ast}$ are three internally disjoint paths of lengths $1,3,3$ between $u_{1}$ and $u^{\ast}$, a contradiction. If $\{u_{2},u_{4}\}\subseteq N_{H^{\ast}}(w)$, then $u^{\ast}u_{2}$,$u_{2}wu_{4}u^{\ast}$,$u_{2}u_{1}u_{3}u^{\ast}$ are three internally disjoint paths of lengths $1,3,3$ between $u^{\ast}$ and $u_{2}$, a contradiction. If $\{u_{4},u_{5}\}\subseteq N_{H^{\ast}}(w)$, then $u^{\ast}u_{4}$,$u_{4}wu_{5}u^{\ast}$,$u_{4}u_{1}u_{3}u^{\ast}$ are three internally disjoint paths of lengths $1,3,3$ between $u^{\ast}$ and $u_{4}$, a contradiction. Thus, $d_{H^{\ast}}(w)\leq1$ for $w\in W_{H^{\ast}}$ and $H^{\ast}\cong S_{r+1}^{1}$ with $r\geq3$. In addition, we can check that $W_{H^{\ast}}\cap W_{N(u^{\ast})\setminus H^{\ast}}=\emptyset$ for $H^{\ast}\cong S_{r+1}^{1}$ with $r\geq2$. Consequently, $d(w)=d_{H^{\ast}}(w)$ for $w\in W_{H^{\ast}}$. Furthermore, let $x_{w_{1}}\geq x_{w_{2}}$, we consider the following two cases.

{\bf Case 1.} At least a vertex $w_{i}\in \cup_{j=1}^{3}N_{W}(u_{j})$ for some $i\in \{1,2\}$.

{\bf Subcase 1.1.} $w_{1}, w_{2}\in \cup_{j=1}^{3}N_{W}(u_{j})$.

In this case, we obtain $(\cup_{j=1}^{3}N_{W}(u_{j}))\cap (W_{N(u^{\ast})}\setminus \{\cup_{j=1}^{3}N_{W}(u_{j})\})=\emptyset$. Thus, $N_{N(u^{\ast})}(w)\subseteq \cup_{j=1}^{3}N_{W}(u_{j})$ for each $w\in \cup_{j=1}^{3}N_{W}(u_{j})$. If $N_{C_{3}}(w_{1})\cap N_{C_{3}}(w_{2})=\emptyset$, then there exists a $\theta(1,3,3)$, a contradiction. If $|N_{C_{3}}(w_{1})\cap N_{C_{3}}(w_{2})|\geq2$, then there exists a $\theta(1,3,3)$, a contradiction. Thus, $|N_{C_{3}}(w_{1})\cap N_{C_{3}}(w_{2})|=1$ and hence there exists a cut vertex, which contradicts Lemma \ref{le:ch-2.4.} {\rm(}$ii${\rm)}.

{\bf Subcase 1.2.} $w_{1}\in \cup_{j=1}^{3}N_{W}(u_{j})$ and $w_{2}\notin \cup_{j=1}^{3}N_{W}(u_{j})$.

In this case, suppose that $w_{2}\in W_{H^{\ast}}$, where $H^{\ast}\cong S_{r+1}^{1}$ with $r\geq3$. Then
\[\begin{cases}
\rho x_{w_{1}}\leq x_{w_{2}}+x_{u_{1}}+x_{u_{2}}+x_{u_{3}}\leq x_{w_{2}}+3x_{u^{\ast}},\\
\rho x_{w_{2}}=x_{w_{1}}+x_{u_{k}}\leq x_{w_{1}}+x_{u^{\ast}},
\end{cases}\]
where $ u_{k}\in V(H^{\ast})\setminus \{u_{1},u_{2},u_{3}\}$. Combining with the above system of inequalities, we get $x_{w_{1}}\leq \frac{3\rho+1}{\rho^{2}-1}x_{u^{\ast}}$. Since $f(x)=\frac{3x+1}{x^{2}-1}$ is decreasing with $x>10$, we have $x_{w_{1}}\leq \frac{3\rho+1}{\rho^{2}-1}x_{u^{\ast}}<\frac{31}{99}x_{u^{\ast}}$. By Lemma \ref{le:ch-2.7.}, we get $1=e(W)<0+\frac{3}{2}-\frac{68}{99}=\frac{161}{198}$, a contradiction. Suppose that $w_{2}\in W_{H^{\star}}$, where $H^{\star}\subseteq N_{+}(u^{\ast})\setminus H^{\ast}$ and $H^{\star}\cong S_{r^{\prime}+1}^{1}$ with $r^{\prime}\geq2$. Furthermore, if $H^{\ast}\cong S_{3}^{1}: u_{1}u_{2}u_{3}$ and $H^{\star}\cong S_{3}^{1}: u_{1}^{\prime}u_{2}^{\prime}u_{3}^{\prime}$, then
\[\begin{cases}
\rho x_{w_{1}}\leq x_{w_{2}}+x_{u_{1}}+x_{u_{2}}+x_{u_{3}}\leq x_{w_{2}}+3x_{u^{\ast}},\\
\rho x_{w_{2}}\leq x_{w_{1}}+x_{u_{1}^{\prime}}+x_{u_{2}^{\prime}}+x_{u_{3}^{\prime}}\leq x_{w_{1}}+3x_{u^{\ast}}.
\end{cases}\]
Combining with the above system of inequalities, we get $ x_{w_{1}}\leq \frac{3}{\rho-1}x_{u^{\ast}}<\frac{1}{3}x_{u^{\ast}}$ due to $\rho>10$. By Lemma \ref{le:ch-2.7.}, we get $1=e(W)<0+\frac{3}{2}-\frac{2}{3}=\frac{5}{6}$, a contradiction. If $H^{\ast}\cong S_{3}^{1}$ and $H^{\star}\cong S_{r^{\prime}+1}^{1}$ with $r^{\prime}\geq3$, then
\[\begin{cases}
\rho x_{w_{1}}\leq x_{w_{2}}+x_{u_{1}}+x_{u_{2}}+x_{u_{3}}\leq x_{w_{2}}+3x_{u^{\ast}},\\
\rho x_{w_{2}}\leq x_{w_{1}}+x_{u^{\ast}}.
\end{cases}
\]
Combining with the above system of inequalities, we get $x_{w_{1}}\leq \frac{3\rho+1}{\rho^{2}-1}x_{u^{\ast}}$. Since $f(x)=\frac{3x+1}{x^{2}-1}$ is decreasing with $x>10$, we have $x_{w_{1}}\leq \frac{3\rho+1}{\rho^{2}-1}x_{u^{\ast}}<\frac{31}{99}x_{u^{\ast}}$. By Lemma \ref{le:ch-2.7.}, we get $1=e(W)<0+\frac{3}{2}-\frac{68}{99}=\frac{161}{198}$, a contradiction. If $H^{\ast}\cong S_{r+1}^{1}$ with $r\geq3$ and $H^{\star} \cong S_{r^{\prime}+1}^{1}$ with $r^{\prime}\geq3$, then
\[\begin{cases}
\rho x_{w_{1}}\leq x_{w_{2}}+x_{u^{\ast}},\\
\rho x_{w_{2}}\leq x_{w_{1}}+x_{u^{\ast}}.
\end{cases}
\]
Combining with the above system of inequalities, we get $x_{w_{1}}\leq \frac{1}{\rho-1}x_{u^{\ast}}<\frac{1}{9}x_{u^{\ast}}$ due to $x>10$. By Lemma \ref{le:ch-2.7.}, we get $1=e(W)<0+\frac{3}{2}-\frac{8}{9}=\frac{11}{18}$, a contradiction. Suppose that $w_{2}\in W_{N_{0}(u^{\ast})}$ and $H^{\ast}\cong S_{3}^{1}$, then
\[\begin{cases}
\rho x_{w_{1}}\leq x_{w_{2}}+x_{u_{1}}+x_{u_{2}}+x_{u_{3}}\leq x_{w_{2}}+3x_{u^{\ast}},\\
\rho x_{w_{2}}\leq x_{w_{1}}+\sum_{u\in N_{0}(u^{\ast})}x_{u}<x_{w_{1}}+\frac{1}{2}x_{u^{\ast}}.
\end{cases}\]
Combining with the above system of inequalities, we get $x_{w_{1}}<\frac{3\rho+\frac{1}{2}}{\rho^{2}-1}x_{u^{\ast}}<\frac{30\frac{1}{2}}{99}$. By Lemma \ref{le:ch-2.7.}, we get $1=e(W)<0+\frac{3}{2}-\frac{68\frac{1}{2}}{99}=\frac{80}{99}$, a contradiction. Suppose that $w_{2}\in W_{N_{0}(u^{\ast})}$ and $H^{\ast}\cong S_{r+1}^{1}$ with $r\geq3$, then
\[\begin{cases}
\rho x_{w_{1}}\leq x_{w_{2}}+x_{u^{\ast}},\\
\rho x_{w_{2}}\leq x_{w_{1}}+\sum_{u\in N_{0}(u^{\ast})}x_{u}<x_{w_{1}}+\frac{1}{2}x_{u^{\ast}}.
\end{cases}\]
Combining with the above system of inequalities, we get $x_{w_{1}}<\frac{2\rho-1}{2\rho^{2}-\rho-2}x_{u^{\ast}}<\frac{19}{188}$. By Lemma \ref{le:ch-2.7.}, we get $1=e(W)<0+\frac{3}{2}-\frac{169}{188}=\frac{113}{188}$, a contradiction.

{\bf Case 2.} $w_{i}\notin \cup_{j=1}^{3}N_{W}(u_{j})$ for any $i\in \{1,2\}$.

Suppose that $w_{1}, w_{2}\in W_{H^{\ast}}$, where $H^{\ast}\cong S_{r+1}^{1}$ for $r\geq3$. Then $N_{H^{\ast}}(w_{1})\cap N_{H^{\ast}}(w_{2})=\emptyset$. Otherwise there exists a cut vertex, which contradicts Lemma \ref{le:ch-2.4.} {\rm(}$ii${\rm)}. By lemma \ref{le:ch-2.4.}, we get $|N_{H^{\ast}}(w_{i})|=1$ for each $i\in \{1,2\}$. Moreover,
\[\begin{cases}
\rho x_{w_{1}}\leq x_{w_{2}}+x_{u^{\ast}},\\
\rho x_{w_{2}}\leq x_{w_{1}}+x_{u^{\ast}}.
\end{cases}\]
Combining with the above system of inequalities, we get $x_{w_{1}}\leq \frac{1}{\rho-1}x_{u^{\ast}}<\frac{1}{9}$ due to $\rho>10$. By Lemma \ref{le:ch-2.7.}, we get $1=e(W)<0+\frac{3}{2}-\frac{8}{9}=\frac{11}{18}$, a contradiction. Suppose that $w_{1}\in W_{H^{\ast}}$, where $H^{\ast}\cong S_{r+1}^{1}$ for $r\geq3$ and $w_{2}\in W_{H^{\star}}$, where $H^{\star}\cong S_{r^{\prime}+1}^{1}$ for $r^{\prime}\geq3$ and $H^{\star}\in N_{+}(u^{\ast})\setminus H^{\ast}$. Then $|N_{H^{\star}}(w_{2})|=1$.
Similar with above, we get a contradiction. Suppose that $w_{1}\in W_{H^{\ast}}$, where $H^{\ast}\cong S_{r+1}^{1}$ for $r\geq3$ and $w_{2}\in W_{N_{0}(u^{\ast})}$. Then $|N_{H^{\ast}}(w_{1})|=1$ and $|N_{N_{0}(u^{\ast})}(w_{2})|\leq |N_{0}(u^{\ast})|$. Thus, we get
\[\begin{cases}
\rho x_{w_{1}}\leq x_{w_{2}}+x_{u^{\ast}},\\
\rho x_{w_{2}}\leq x_{w_{1}}+\sum_{u\in N_{0}(u^{\ast})}x_{u}<x_{w_{1}}+\frac{1}{2}x_{u^{\ast}}.
\end{cases}\]
Combining with the above system of inequalities, we get $x_{w_{1}}\leq \frac{\rho+\frac{1}{2}}{\rho^{2}-1}x_{u^{\ast}}<\frac{10\frac{1}{2}}{99}$. By Lemma \ref{le:ch-2.7.}, we get $1=e(W)<0+\frac{3}{2}-\frac{88\frac{1}{2}}{99}=\frac{60}{99}$, a contradiction. Suppose that $w_{1}, w_{2}\in W_{N_{0}(u^{\ast})}$ for each $i\in \{1,2\}$. Then $|N_{N_{0}(u^{\ast})}(w_{i})|\leq |N_{0}(u^{\ast})|$ for each $i\in \{1,2\}$. Thus, we get
\[\begin{cases}
\rho x_{w_{1}}\leq x_{w_{2}}+\sum_{u\in N_{0}(u^{\ast})}x_{u}<x_{w_{1}}+\frac{1}{2}x_{u^{\ast}},\\
\rho x_{w_{2}}\leq x_{w_{1}}+\sum_{u\in N_{0}(u^{\ast})}x_{u}<x_{w_{1}}+\frac{1}{2}x_{u^{\ast}}.
\end{cases}\]
Combining with the above system of inequalities, we get $x_{w_{1}}<\frac{1}{2(\rho-1)}x_{u^{\ast}}<\frac{1}{18}x_{u^{\ast}}$. By Lemma \ref{le:ch-2.7.}, we get $1=e(W)<0+\frac{3}{2}-\frac{17}{18}=\frac{5}{9}$, a contradiction. This completes the proof of Lemma \ref{le:ch-3.4.}. $\qedsymbol$

\noindent\begin{lemma}\label{le:ch-3.5.} $G^{\ast}[N(u^{\ast})]$ contains no triangle.
\end{lemma}
\noindent{\bf Proof.} Suppose on the contrary that $G^{\ast}[N(u^{\ast})]$ contains triangles. Then $G^{\ast}[N(u^{\ast})]$ contains a component which is isomorphic to $S_{r+1}^{1}$ for $r\geq2$. Let $H^{\ast}\cong S_{r+1}^{1}$ be a component of $G^{\ast}[N(u^{\ast})]$. Then $e(H^{\ast})=r+1$. Let $u_{1}u_{2}u_{3}u_{1}$ be the triangle of $H^{\ast}$ and $d_{H^{\ast}}(u_{1})=d_{H^{\ast}}(u_{2})$. If $W_{H^{\ast}}=\emptyset$, then $x_{u_{1}}=x_{u_{2}}$. Furthermore,
$$\rho x_{u_{1}}=x_{u_{2}}+x_{u_{3}}+x_{u^{\ast}}\leq x_{u_{1}}+2x_{u^{\ast}}.$$ Thus, $x_{u_{1}}\leq \frac{2}{\rho-1}x_{u^{\ast}}<\frac{2}{9}x_{u^{\ast}}$ due to $\rho>10$.
$$
\begin{aligned}
\zeta(H^{\ast})&=x_{u_{1}}+x_{u_{2}}+(r-1)x_{u_{3}}
\\&<(\frac{4}{9}+r-1)x_{u^{\ast}}
\\&=(e(H^{\ast})-\frac{14}{9})x_{u^{\ast}}.
\end{aligned}
$$
Recall that $\zeta(H)\leq e(H)x_{u^{\ast}}$ for $H\in \mathcal{H}\setminus H^{\ast}$. Thus, $$
\begin{aligned}
\sum_{H\in \mathcal{H}}\zeta(H)&=\zeta(H^{\ast})+\sum_{H\in \mathcal{H}\setminus H^{\ast}}\zeta(H)
\\&<(e(H^{\ast})-\frac{14}{9})x_{u^{\ast}}+\sum_{H\in \mathcal{H}\setminus H^{\ast}}e(H)x_{u^{\ast}}
\\&=(e(N_{+}(u^{\ast}))-\frac{14}{9})x_{u^{\ast}}
\end{aligned}
,$$
which contradicts \eqref{eq:ch-3}. Thus, $W_{H^{\ast}}\neq\emptyset$.

Since $e(W)=0$, combining Lemma \ref{le:ch-2.4.}, we have $d_{N(u^{\ast})}(w)\geq2$. If $r\geq3$, then $W_{H^{\ast}}\cap W_{N(u^{\ast})\setminus H^{\ast}}=\emptyset$, and hence $d(w)=d_{H^{\ast}}(w)\leq1$ for $w\in W_{H^{\ast}}$, a contradiction. Thus, $r=2$, that is $H^{\ast}$ is a triangle $u_{1}u_{2}u_{3}$. Since $G^{\ast}$ is $\theta(1,3,3)$-free, we obtain that $W_{H^{\ast}}\cap W_{N(u^{\ast})\setminus H^{\ast}}=\emptyset$, and hence $d(w)=d_{H^{\ast}}(w)\leq3$. First, we assume that $|W_{H^{\ast}}|=1$. Let $W_{H^{\ast}}=\{w\}$. As $H^{\ast}$ is a triangle, we obtain $2\leq d(w)=d_{H^{\ast}}(w)\leq3$. We consider two cases as follows.

{\bf Case 1.} $d_{H^{\ast}}(w)=2$.

In this case, without loss of generality, we suppose $N(w)=\{u_{1},u_{2}\}$. Then $x_{u_{1}}=x_{u_{2}}$. Since
$$\rho x_{u_{3}}=x_{u_{1}}+x_{u_{2}}+x_{u^{\ast}}\leq 3x_{u^{\ast}}.$$
Furthermore, $x_{u_{3}}\leq \frac{3}{\rho}x_{u^{\ast}}$.
Similarly, $$\rho x_{u_{1}}=x_{u_{2}}+x_{u_{3}}+x_{u^{\ast}}+x_{w}\leq x_{u_{1}}+\frac{3}{\rho}x_{u^{\ast}}+2x_{u^{\ast}},$$
we obtain that $x_{u_{1}}\leq \frac{2\rho+3}{\rho(\rho-1)}x_{u^{\ast}}$. Thus, $$\zeta(H^{\ast})=x_{u_{1}}+x_{u_{2}}+x_{u_{3}}\leq \frac{7\rho+3}{\rho(\rho-1)}x_{u^{\ast}}.$$ Since
$\frac{7x+3}{x(x-1)}$ is decreasing in variable $x>10$, we get
$$\zeta(H^{\ast})<\frac{73}{90}x_{u^{\ast}}<(e(H^{\ast})-\frac{3}{2})x_{u^{\ast}}.$$
Recall that $\zeta(H)\leq e(H)x_{u^{\ast}}$ for $H\in \mathcal{H}\setminus H^{\ast}$. Thus,
$$
\begin{aligned}
\sum_{H\in \mathcal{H}}\zeta(H)&=\zeta(H^{\ast})+\sum_{H\in \mathcal{H}\setminus H^{\ast}}\zeta(H)
\\&<(e(H^{\ast})-\frac{3}{2})x_{u^{\ast}}+\sum_{H\in \mathcal{H}\setminus H^{\ast}}\zeta(H)
\\&=(e(N_{+}(u^{\ast}))-\frac{3}{2})x_{u^{\ast}}
\end{aligned}
,$$
which contradicts \eqref{eq:ch-3}.

{\bf Case 2.} $d_{H^{\ast}}(w)=3$.

In this case, let $N_{H^{\ast}}(w)=\{u_{1},u_{2},u_{3}\}$. Then $x_{u_{1}}=x_{u_{2}}=x_{u_{3}}$. Since
$$\rho x_{u_{1}}=x_{u_{2}}+x_{u_{3}}+x_{u^{\ast}}+x_{w}\leq 2x_{u_{1}}+2x_{u^{\ast}},$$ we obtain $x_{u_{1}}\leq \frac{2}{\rho-2}x_{u^{\ast}}$. Furthermore,
$$\zeta(H^{\ast})=x_{u_{1}}+x_{u_{2}}+x_{u_{3}}\leq \frac{6}{\rho-2}x_{u^{\ast}}<\frac{3}{4}x_{u^{\ast}}<(e(H^{\ast})-\frac{3}{2})x_{u^{\ast}}.$$ Thus,
$$\sum_{H\in \mathcal{H}}\zeta(H)<(e(N_{+}(u^{\ast}))-\frac{3}{2})x_{u^{\ast}},$$ a contradiction. Thus, $|W_{H^{\ast}}|\geq2$.

Recall that $W_{H^{\ast}}\cap W_{N(u^{\ast})\setminus H^{\ast}}(w)=\emptyset$ for any $w\in W_{H^{\ast}}$ and $2\leq d(w)=d_{H^{\ast}}(w)\leq3$. Let $w_{1}\in W_{H^{\ast}}$ such that $N(w_{1})=\{u_{1},u_{2},u_{3}\}$ and $w_{2}\neq w_{1}$ in $W_{H^{\ast}}$ satisfying $d(w_{2})\geq2$. Suppose that $u_{1},u_{2}\in N(w_{2})$. Then $u^{\ast}u_{1}, u^{\ast}u_{3}w_{1}u_{1}$ and $u^{\ast}u_{2}w_{2}u_{1}$ are three internally disjoint paths of lengths $1,3,3$ between $u^{\ast}$ and $u_{1}$, a contradiction. Hence, $d(w)=d_{H^{\ast}}(w)=2$ for any $w\in W_{H^{\ast}}$. This implies $1\leq |N(w_{1})\cap N(w_{2})|\leq2$ for any two vertices $w_{1},w_{2}\in W_{H^{\ast}}$. If $|N(w_{1})\cap N(w_{2})|=1$, without loss of generality, let $N(w_{1})=\{u_{1},u_{2}\}$ and $N(w_{2})=\{u_{1},u_{3}\}$, then $u^{\ast}u_{1}, u^{\ast}u_{3}w_{2}u_{1}$ and $u^{\ast}u_{2}w_{1}u_{1}$ are three internally disjoint paths of lengths $1,3,3$ between $u^{\ast}$ and $u_{1}$, a contradiction. Thus, $|N(w_{1})\cap N(w_{2})|=2$. That is $N(w_{1})=N(w_{2})$ for any two vertices $w_{1},w_{2}\in W_{H^{\ast}}$. Without loss of generality, suppose that $N(w)=\{u_{1},u_{2}\}$ for any $w\in W_{H^{\ast}}$. Let $G_{1}=G^{\ast}-\{u_{1}w| w\in N_{W}(u_{1})\}+\{u^{\ast}w| w\in N_{W}(u_{1})\}$. Obviously, $G_{1}$ is $\theta(1,3,3)$-free. By Lemma \ref{le:ch-2.1.}, we get $\rho(G_{1})>\rho$, a contradiction. This completes the proof of Lemma \ref{le:ch-3.5.}. $\qedsymbol$

\noindent{\bf Proof of Theorem 1.4.} By Lemmas \ref{le:ch-3.1.},\ref{le:ch-3.3.} and \ref{le:ch-3.5.}, we obtain that each component of $G^{\ast}[N_{+}(u^{\ast})]$ is a non-trivial tree. If $c=0$, then $G^{\ast}$ is bipartite. By Lemma \ref{le:ch-2.5.}, we have $\rho\leq\sqrt{m}<\frac{1+\sqrt{4m-3}}{2}$ for $m\geq 92$, a contradiction. Thus $c=1$ and $\sum_{u\in N_{0}(u^{\ast})}\frac{x_{u}}{x_{u^{\ast}}}<\frac{1}{2}$ from Inequality \eqref{eq:ch-4}. Let $H$ be the unique component of $G^{\ast}[N_{+}(u^{\ast})]$, where $H$ is a non-trivial tree. By Lemma \ref{le:ch-3.1.} {\rm(}$iii${\rm)} and {\rm(}$iv${\rm)}, we have $diam(H)\leq3$.

If $diam(H)=3$, then $H$ is a double star. Let $u_{1}$ and $u_{2}$ be the two center vertices of $H$. Let $\{v_{1}, v_{2}, \cdots, v_{a}\}\in N_{H}(u_{1})\setminus u_{2}$ and let $\{z_{1}, z_{2}, \cdots, z_{b}\}\in N_{H}(u_{2})\setminus u_{1}$ for $a, b\geq1$. If $W_{H}=\emptyset$, without loss of generality, assume that $x_{u_{1}}\geq x_{u_{2}}$, then let $G_{2}=G^{\ast}-\{u_{2}v| v\in N_{H}(u_{2})\setminus \{u_{1}\}\}+\{u_{1}v| v\in N_{H}(u_{2})\setminus \{u_{2}\}\}$. We can verify that $G_{2}$ is $\theta(1,3,3)$-free. By Lemma \ref{le:ch-2.1.}, we get $\rho(G_{2})> \rho$, a contradiction. Thus, $W_{H}\neq\emptyset$. It is easily checked $W_{H}\cap W_{N(u^{\ast})\setminus H}=\emptyset$ for any $w\in W_{H}$. Hence, $N(w)\subseteq V(H)$ and by Lemma \ref{le:ch-2.4.}, $d(w)=d_{H}(w)\geq2$. We claim $d(w)=d_{H}(w)=2$. Suppose contrary that $d_{H}(w)\geq3$. If $\{u_{1},u_{2},v_{1}\}\in W_{H}(w)$, then $u_{1}u_{2}$, $u_{1}v_{1}wu_{2}$, $u_{1}u^{\ast}z_{1}u_{2}$ are three internally disjoint paths of lengths $1,3,3$ between $u_{1}$ and $u_{2}$, a contradiction. If $\{u_{1},v_{1},v_{2}\}\in W_{H}(w)$, then $u^{\ast}v_{2}$, $u^{\ast}v_{1}wv_{2}$, $u^{\ast}u_{2}u_{1}v_{2}$ are three internally disjoint paths of lengths $1,3,3$ between $u^{\ast}$ and $v_{2}$, a contradiction. If $\{u_{2},v_{1},v_{2}\}\in W_{H}(w)$, then $u^{\ast}u_{2}$, $u^{\ast}v_{1}wu_{2}$, $u^{\ast}v_{2}u_{1}u_{2}$ are three internally disjoint paths of lengths $1,3,3$ between $u^{\ast}$ and $u_{2}$, a contradiction. If $\{v_{1},v_{2},z_{1}\}\in W_{H}(w)$, then $v_{1}u^{\ast}$, $v_{1}wv_{2}u^{\ast}$, $v_{1}u_{1}u_{2}u^{\ast}$ are three internally disjoint paths of lengths $1,3,3$ between $v_{1}$ and $u^{\ast}$, a contradiction. If $\{v_{1},v_{2}
,v_{3}\}\in W_{H}(w)$, then $u^{\ast}v_{1}$, $u^{\ast}v_{3}wv_{1}$, $u^{\ast}u_{2}u_{1}v_{1}$ are three internally disjoint paths of lengths $1,3,3$ between $u^{\ast}$ and $v_{1}$, a contradiction. Thus, $d(w)=d_{H}(w)=2$. We claim $N(w)=\{u_{1},u_{2}\}$. Otherwise, we consider five cases in the sense of symmetry as follows. If $N(w)=\{u_{1}, v_{1}\}$, then $u^{\ast}u_{1}$, $u^{\ast}v_{1}wu_{1}$, $u^{\ast}z_{1}u_{2}u_{1}$ are three internally disjoint paths of lengths $1,3,3$ between $u^{\ast}$ and $u_{1}$, a contradiction. If $N(w)=\{u_{1}, z_{1}\}$, then $u_{1}u_{2}$, $u_{1}wu^{\ast}z_{1}u_{2}$, $u_{2}u^{\ast}z_{1}u_{1}$ are three internally disjoint paths of lengths $1,3,3$ between $u_{1}$ and $u_{2}$, a contradiction. If $N(w)=\{v_{1}, v_{2}\}$, then $u^{\ast}v_{1}$, $u^{\ast}v_{2}wv_{1}$, $u^{\ast}u_{2}u_{1}v_{1}$ are three internally disjoint paths of lengths $1,3,3$ between $u^{\ast}$ and $v_{1}$, a contradiction. If $N(w)=\{v_{1}, z_{1}\}$, then $u^{\ast}v_{1}$, $u^{\ast}z_{1}wv_{1}$, $u^{\ast}u_{2}u_{1}v_{1}$ are three internally disjoint paths of lengths $1,3,3$ between $u^{\ast}$ and $v_{1}$, a contradiction. Thus, $N(w)=\{u_{1},u_{2}\}$. Let $G_{3}=G^{\ast}-\{u_{2}w| w\in N_{W}(u_{2})\}+\{u^{\ast}w| w\in N_{W}(u_{2})\}$. We can verify that $G_{3}$ is $\theta(1,3,3)$-free. By Lemma \ref{le:ch-2.1.}, we get $\rho(G_{3})>\rho$, a contradiction.

If $diam(H)\leq2$, then $H\cong K_{1,r}$ with $r\geq1$. Let $V(H)=\{u_{0},u_{1},\cdots, u_{r}\}$ and $u_{0}$ be the center vertex of $H$ with $r\geq1$. We claim $r\geq 9$. By $\rho>10$ and Inequality \eqref{eq:ch-4}, we have $$10x_{u^{\ast}}<\rho x_{u^{\ast}}=x_{u_{0}}+x_{u_{1}}+\cdots+x_{u_{r}}+\sum_{v\in N_{0}(u^{\ast})}x_{v} < (r+1+\frac{1}{2})x_{u^{\ast}}.$$ Thus, $r\geq9$. We claim $W_{H}=\emptyset$. Suppose on the contrary that $W_{H}\neq\emptyset$.

First, we assume $d_{H}(w)\geq 3$ for any vertex $w\in W_{H}$. If $\{u_{0}, u_{1},u_{2}\}\in N_{H}(w)$, then $u^{\ast}u_{1}, u^{\ast}u_{2}wu_{1}$ and $u^{\ast}u_{3}u_{0}u_{1}$ are three internally disjoint paths of lengths $1,3,3$ between $u^{\ast}$ and $u_{1}$, a contradiction. If $\{u_{1},u_{2}, u_{3}\}\in N_{H}(w)$, then $u^{\ast}u_{3}, u^{\ast}u_{4}u_{0}u_{3}$ and $u^{\ast}u_{1}wu_{3}$ are three internally disjoint paths of lengths $1,3,3$ between $u^{\ast}$ and $u_{3}$, a contradiction. Thus, $d_{H}(w)\leq 2$. If $d_{H}(w)=1$, then we obtain that $w$ is only adjacent to the center vertex $u_{0}$. Otherwise, let $N_{H}(w)={u_{1}}$. By Lemma \ref{le:ch-2.4.} and $e(W)=0$, we obtain $|N_{N_{0}(u^{\ast})}(w)|\geq1$. Then $u^{\ast}u_{1}$, $u^{\ast}vwu_{1}$, $u^{\ast}u_{2}u_{0}u_{1}$ are three internally disjoint paths of lengths $1,3,3$ between $u^{\ast}$ and $u_{1}$, where $v\in N_{N_{0}(u^{\ast})}(w)$, a contradiction. Thus, $N_{H}(w)=\{u_{0}\}$. By Inequality \eqref{eq:ch-4}, we have
$$\rho x_{w}\leq x_{u_{0}}+\sum_{v\in N_{0}(u^{\ast})}x_{v}<\frac{3}{2}x_{u^{\ast}}.$$
Thus, $x_{w}<\frac{3}{2\rho}x_{u^{\ast}}<\frac{3}{20}x_{u^{\ast}}$ due to $\rho> 10$. By Lemma \ref{le:ch-2.7.}, we get $0=e(W)<-1+\frac{3}{2}-\frac{17}{20}=-\frac{7}{20}$, a contradiction. If $d_{H}(w)=2$, then we have $N_{N_{0}(u^{\ast})}(w)=\emptyset$. Otherwise, $G^{\ast}$ contains a $\theta(1,3,3)$. Furthermore, we get
$$\rho x_{w}\leq x_{u_{0}}+x_{u_{1}}\leq 2x_{u^{\ast}}.$$ This implies that $x_{w}\leq \frac{1}{2\rho}x_{u^{\ast}}< \frac{1}{20}x_{u^{\ast}}$ due to $\rho> 10$. By Lemma \ref{le:ch-2.7.}, we get $0=e(W)<-1+\frac{3}{2}-\frac{19}{20}<0$, a contradiction. Thus, $W_{H}=\emptyset$. If $W\neq \emptyset$, then by Lemma \ref{le:ch-2.1.}, we obtain $d(w)=d_{N_{0}(u^{\ast})}(w)$ for any vertex $w\in W$. Furthermore, $N_{0}(u^{\ast})\neq \emptyset$ and $$\rho x_{w}\leq \sum_{v\in N_{0}(u^{\ast})}x_{v}< \frac{1}{2}x_{u^{\ast}}.$$ This implies that
$x_{w}<\frac{1}{2\rho}x_{u^{\ast}}<\frac{1}{20}x_{u^{\ast}}$ due to $\rho> 10$. By Lemma \ref{le:ch-2.7.}, we get $0=e(W)< -1+\frac{3}{2}-\frac{19}{20}=-\frac{9}{20}$, a contradiction. Thus, $W=\emptyset$ and hence $G^{\ast}\cong G_{4}$ (see Figure. 1). Let $|N_{0}(u^{\ast})|=t$. Since $m$ is even, we obtain that $t$ is odd and $t\geq1$. By Lemma \ref{le:ch-2.3.}, we obtain that $\rho$ is the largest root of the equation $f(x,t)=0$ where $$f(x,t)=x^{4}-mx^{2}-(m-t-1)x+\frac{t(m-t-1)}{2}$$ for $m=t+1+2r\geq92$. Since $$f(x,t)-f(x,1)=(t-1)x+\frac{(t-1)(m-t-2)}{2}>0$$ for $x>0$ and $t\geq3$, which implies that $t=1$ for the extremal graph $G^{\ast}$. By Lemma \ref{le:ch-2.6.}, we have $\rho(S_{\frac{m+4}{2},2}^{-})> \frac{1+\sqrt{4m-5}}{2}$ for $ m\geq92$ and $G^{\ast}\cong S_{\frac{m+4}{2},2}^{-}$, as desired. This completes the proof of Theorem 1.4. $\blacksquare$
\begin{figure}[H]
\begin{centering}
\includegraphics[scale=0.25]{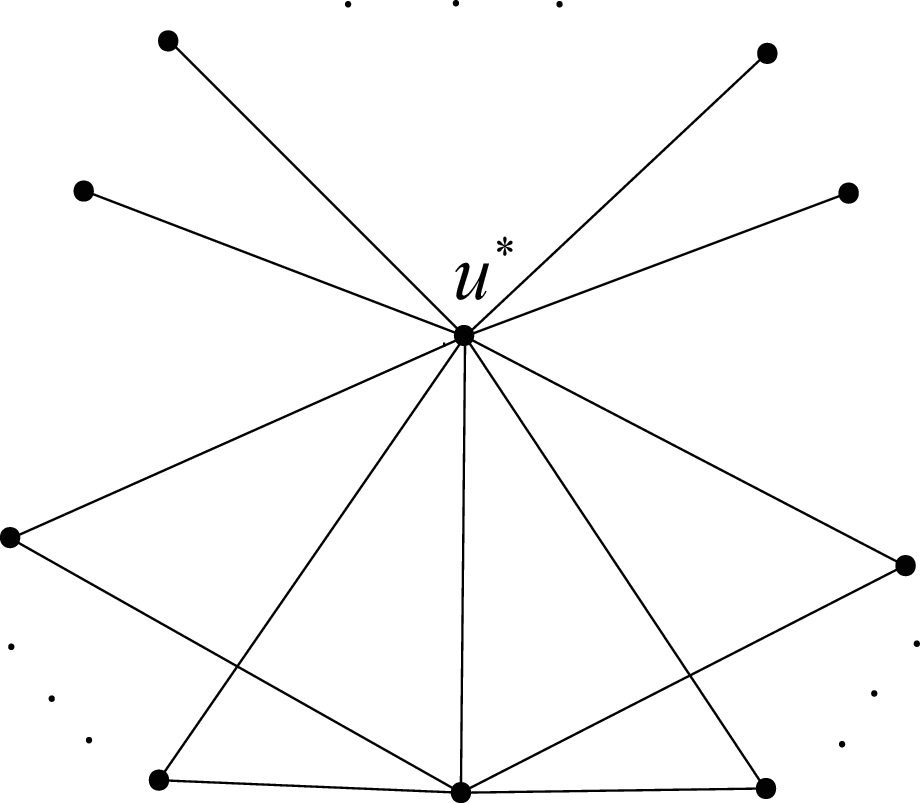}
\caption{The graph $G_{4}$.}\label{fi:ch-1}
\end{centering}
\end{figure}

\section*{Data availability}

No data was used for the research described in the article.

\section*{Declaration of competing interest}

The authors declare that they have no conflict of interest.


\begin{thebibliography}{99}
\bibitem{BoMu1} J.A. Bondy, U.S.R. Murty, Graph Theory, in: Graduate Texts in Mathematics, Vol. 244, Springer, New York, 2008.
\bibitem{BrH} R.A. Brualdi, A.J. Hoffman, On the spectral radius of $(0,1)$ matrices, Linear Algebra Appl. 65 (1985)
133--146.
\bibitem{ChZ} M.Z. Chen, X.-D Zhang, Some new results and problems in spectral extremal graph theory, J. Anhui Univ. Nat. Sci. 42 (2018) 12--25 (in Chinese).
\bibitem{CvRS} D. Cvektovi\'{c}, P. Rowlinson, S. Simi\'{c}, An Introduction to the Theory of Graph Spectra, Cambrige University Press, New York, 2010.
\bibitem{FaY} X.N. Fang, L.H. You, The maximum spectral radius of graphs of given size with forbidden subgraph, Linear Algebra Appl. 666 (2023) 114--128.
\bibitem{FSi} Z. F\"{u}redi, M. Simonovits, The history of degenerate (bipartite) extremal graph problems, Erd\H{o}s centennial, Bolyai Soc. Math. Stud. 25 (2013) 169--264.
\bibitem{GaoL} J. Gao, X.L. Li, The maximum spectral radius of $\theta(1,3,3)$-free graphs with given size, arXiv: 2410.07721v1.
\bibitem{LiZZ} S.C. Li, S.S. Zhao, L.T. Zou, Spectral extrema of graphs with fixed size: forbidden fan graph,
friendship graph or theta graph, arXiv: 2409.15918v1.
\bibitem{LiZS} X. Li, M. Zhai, J. Shu, A Brualdi-Hoffman-Tur\'{a}n problem on cycles, European J. Combin. 120 (2024) 103966.
\bibitem{LiLF} Y.T. Li, W.J. Liu, L.H. Feng, A survey on spectral conditions for some extremal graph problems, Adv. Math. (China) 51 (2022) 193--258.
\bibitem{LuLL} J. Lu, L. Lu, Y. Li, Spectral radius of graphs forbidden $C_{7}$ or $C_{6}^{\Delta}$, Discrete Math. 347(2024) 113781.
\bibitem{LiuW} Y.X. Liu, L.G. Wang, Spectral radius of graphs of given size with forbidden subgraphs, Linear Algebra Appl. 689 (2024) 108--125.
\bibitem{MiLH} G. Min, Z.Z. Lou, Q.X. Huang, A sharp upper bound on the spectral radius of $C_{5}$-free/$C_{6}$-free graphs with given size, Linear Algebra Appl. 640 (2022) 162--178.
\bibitem{Ni5} V. Nikiforov, Some inequalities for the largest eigenvalue of a graph, Combin. Probab. Comput. 11 (2002) 179--189.
\bibitem{Ni6} V. Nikiforov, Walks and the spectral radius of graphs, Linear Algebra Appl. 418 (2006) 257--268.
\bibitem{Ni2} V. Nikiforov, The spectral radius of graphs without paths and cycles of specified length, Linear Algebra Appl. 432 (2010) 2243--2256.
\bibitem{Ni3} V. Nikiforov, Some new results in extremal graph theory. Lond. Math. Soc. Lect. Note Ser. 392 (2011) 141--181.
\bibitem{Nos} E. Nosal, Eigenvalues of Graphs, Master's Thesis, University of Calgary, 1970.
\bibitem{SLW} W.T. Sun, S.C. Li, W. Wei, Extensions on spectral extrema of $C_{5}/C_{6}$-free graphs with given size, Discrete Math. 346 (2023) 113591.
\bibitem{ZhLS} M.Q. Zhai, H.Q. Lin, J.L. Shu, Spectral extrema of graphs with fixed size: Cycles and complete bipartite graphs, European J. Combin. 95 (2021) 103322.
\bibitem{ZhWF1} M.Q. Zhai, B. Wang, Proof of a conjecture on the spectral radius of $C_{4}$-free graphs, Linear Algebra Appl. 437 (2012) 1641--1647.
\end{thebibliography}
\end{document}